# Layer Groups of Batak Weaves

Ma. Louise Antonette De Las Peñas[a], Agnes Garciano[a],
Debbie Marie Verzosa[b] and Mark Tomenes[a]

[a]Department of Mathematics, Ateneo de Manila University, Loyola Heights, Quezon City, Metro Manila, 1108, Philippines
[b]Department of Mathematics and Statistics, University of Southern Mindanao, Kabacan, Cotabato, 9407, Philippines

**Abstract.** This paper discusses the symmetry structures of the decorative weave patterns in baskets and other household items of the Batak, an indigenous community in the Philippines. The study confirms the realization of 15 layer groups that occur as symmetry groups of the Batak weaves. Layers groups are crystallographic space groups that have translational symmetries in two dimensions.

## I. Introduction

Basketry has a long cultural history, spanning thousands of years (Gilhus, 2020). Baskets from various cultures show unique patterns and designs, making them intriguing subjects for mathematical symmetry investigations.

Within the Philippines, diverse cultures developed distinct basket weaving traditions, including the *apugan* of the Kabihug in Camarines Norte (Rubio, 2016), the *bakat* of the Sugbuanon in Cebu and Bohol (Inocian, *et al.*, 2019), and the *Tingkep* of the Pala'wan in Palawan (Dressler, 2019). There is also the long basket weaving tradition of the Batak. This community is spread out in settlements in the island of Palawan (City Government of Puerto Princesa, 2009). They create baskets which, aside from being functional, are also pieces of art, containing elaborate black and white patterns that demonstrate the weavers' technical expertise and understanding of algebraic and geometric principles (De las Peñas, *et al.*, 2021). There are detailed records of Batak woven patterns (Calderon, 1986; Novellino, 2009) and ethnographic accounts of Batak basket-weaving knowledge (Novellino, 2009).

In this paper, we investigate the presence of layer group symmetries in the weave patterns woven by the Batak, which appear predominantly in baskets. In particular, we determine which of the 80 layer group structures known in Crystallography (Kopský & Litvin, 2002) are present in the Batak weave patterns. This work continues from a previous study (De Las Peñas *et al.*, 2021) where a detailed discussion of the weaving techniques of the Batak was presented, including how the decorative repeating patterns in the Batak weaves were created from these techniques.

The Batak patterns analyzed in this paper were culled from existing records (Calderon, 1986; Novellino, 2009), museums, and exhibits (the Batak Visitor Center, the National Commission for Culture and the Arts, the Cotabato museum). Basket designs were also obtained from social enterprises focusing on Batak baskets, and from a two-day visit to a Batak community in Palawan.

This article contributes to existing works on symmetries of culture, particularly on the applications of crystallographic groups to study symmetries of artwork of indigenous communities. See for example the works by Washburn and Crowe and references therein (Washburn & Crowe, 1988, 2004), and other related studies on symmetries of fabrics (De Las Peñas & Amores, 2016; De Las Peñas *et al.*, 2018; Pebryani, 2018; Chudasri & Sukantamala, 2023). Moreover, symmetry analyses have been performed on basket designs of the Yurok, Karok and Hupa indigenous people of California (Washburn, 1986), the Tlingit people of Alaska and Canada (Mukhopadhyay, 2009), and those found in the United States Southwest (Washburn & Webster, 2006). The work described in the present paper, which investigates the presence of layer group structures in Batak baskets, adds to the aforementioned existing studies on symmetries of baskets. Layer groups have been discussed as symmetries of woven fabrics in the literature (Roth, 1993; De Las Peñas *et al.*, 2023), but not in the context of cultural artwork such as baskets. This paper also intends to contribute to studies in the mathematics of baskets. The ethno-mathematician Gerdes discussed the mathematics that emerges in the basketry of Mozambique and of the Peruvian and Brazilian Amazon (Gerdes, 2004, 2010, 2011).

The organization of this paper is as follows. Section 2 presents the Batak weaving process. Section 3 gives the geometric setting with regards to weaves, and explains how to carry out the analysis of their symmetry structures. Section 4 presents the layer groups found in the Batak weaves. Finally, Section 5 gives the concluding remarks and future direction of the work.

**II. Weaving Process**

Batak baskets are made from a variety of bamboo that grows abundantly in the mountainous villages of Puerto Princesa, where Batak communities are found. Two sets of thin bamboo strips are used in basket-making to arrive at black and white patterns: a set of bamboo strips in their natural color ("white-like"); and a set of one-sided black bamboo strips from a blackened pole, obtained by holding the bamboo pole over a flame of a lighted resin, until soot covers its outer portion. Details of the preparation of these white and black bamboo strips as well as the weaving process itself which will be described below, may be found in the literature (Calderon 1986; Novellino, 2009).

Weaving is predominantly done by women. A Batak weaver begins by weaving the base of a basket. First, a set of thin, white bamboo strips are laid out vertically on the floor. The weaver then inserts white strips over and under the vertical ones, interlacing them perpendicularly so that these do not fall apart. The plaited white strands will form the base of the basket. To weave the sides of the basket, the weaver uses a finished woven basket to serve as a model or template for the new one (Fig. 1(a)). The model basket is turned upside down and the white plaited strands are then put on top of the base of the model basket and tied to it temporarily. Then four sticks are secured and tied on each corner of the model basket to form the frame of the sides of the new basket (Novelino, 2009).

In order to start weaving the sides of the new basket, the strands hanging from the base are bent vertically towards the sides of the finished basket (Fig.1(b)). Thus, from a two-dimensional planar weave, the basket weaver now transitions into working out a three-dimensional craft. All the strands used in making the base of a basket now become the vertical strands of the sides of the basket. To make a basket with a black and white decorative pattern that appears on the basket's sides, one-sided black strands are inserted as horizontal strands, with the black part of the strands facing front. In weaving, the vertical and horizontal strands are called *warp* and *weft* strands, respectively. The weft strands are woven continuously around the sides of the basket until a basket of a desired height is achieved. Finally, the four sticks are removed and rattan strands are entwined using a cross stitch pattern around a space, approximately one inch below the rim as well as around the rim of the basket itself in order to secure the basket's opening.

A characteristic of Batak baskets is that the base is partitioned into quadrants; hence the base consists of an even number of vertical and horizontal strands. Two perpendicular lines (quadrantal axes) resulting from the weaving process, are apparent at the base (Fig. 1 (c)). When the base strands are folded up to form the sides of a basket, these quadrantal axes extend and transform into vertical lines through the middle on each of the four sides of the basket. Because of the nature of the Batak weaving technique, warp strands are always colored white while the weft strands are always black. Furthermore, at any step of the weaving exercise, the warp and weft strands are perpendicular, thus either a black portion of a strand covers a white strand or vice versa. The weaver precisely determines the number of strands in each quadrant of the base, based on the specific black and white pattern she intends to create. (De Las Peñas *et al.*, 2020). For example, to arrive at the *Piyakdan* pattern in the basket shown in Fig. 1(d), the weaver uses $16 \times 16$ strands per quadrant (Fig.1(c)).

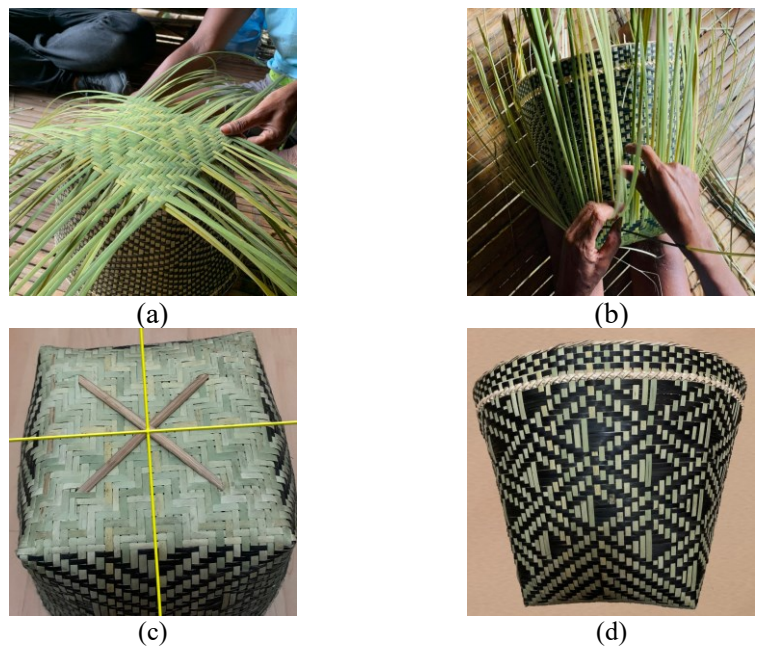

(a) (b)
(c) (d)

Figure 1: (a) Plaited white strands temporarily tied on the base of the model basket; (b)strands in the base being bent towards the sides of the model basket then black horizontal strands inserted to create the black and white pattern; (c) base of the basket showing quadrantal axes; and (d) basket showing the *Piyakdan* pattern.

Although the Batak are known for their baskets, the Batak also weave non-baskets such as trays, and mats. For the non-baskets, the Batak use one-sided black strips. On one side, all strands of the same direction will have the same color. For example, all warp stands at the front and back of the non-basket will be black, and white, respectively. Then all weft strands at the front and back will be white and black, respectively. This way, the same decorative patterns appear at the front and back (see for example, Figs. 2(a)-(b)).

One of the weaving techniques used by the Batak to arrive at the repeating black and white patterns in non-baskets, can be likened to the concept of *thick stripping* (Roth, 1995) where a fixed number $x$ of adjacent weft strands alternate colors black and white; and the same goes for the warp strands. Consider for example two variants of the *Binalang* pattern shown in Figs. 2(c)-(d). In these trays, the number $x$ is changed to show different variants of a pattern, $x$ takes on the values 3 (Fig. 2(c)) and 6 (Fig. 2(d)) respectively. The underlying weave (without colors) is a twill – in particular, a 2-over and 2-under weave of warp over weft (and vice versa).

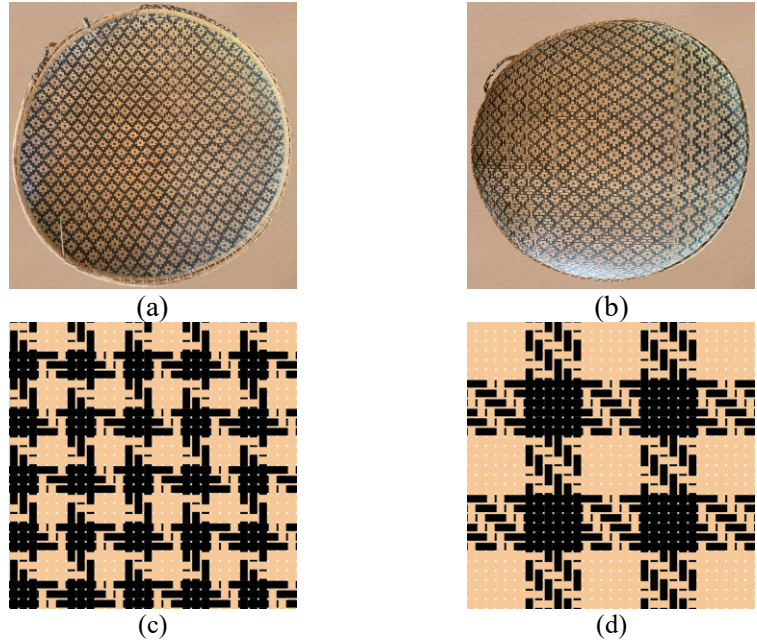

Figure 2: (a)-(b) Tray showing the same decorative pattern in front and back; kind of twill using thick stripping of (c) three white and three black strands; and (d) six white and six black strands alternating in the warp and the weft of two Batak mats.

### III. Mathematical Considerations

There are several weaving techniques employed in the creation of baskets such as the 2-way, 2-fold weave; the 3-way, 3-fold weave; and the *kagome* (hexagon-triangle-hexagon) pattern (Tarnai, 2006). The Batak use the 2-way and 2-fold weave in constructing their baskets. In a *2-way weave*, the strands have two directions: warp (horizontal) and weft (vertical); and are at right angles to each other. A *2-fold weave* means there are no more than two strands crossing each other (Hoskins & Thomas, 1991).

One way to understand the geometric properties of a 2-way 2-fold weave $W$ is to represent it by means of its *design*, $D(W)$. Assuming that the weave repeats indefinitely in the plane via two linearly independent translations, and is periodic, we use the regular tiling of the plane by unit squares. Each square is the intersection of a row of squares (a warp strand) and a column of squares (a weft strand). The square is given a particular color (in this paper, grey) if the weft strand passes over the warp strand at the given intersection, and another color (white) if the warp strand passes over the weft strand. Any weave that we assume to be periodic, therefore, can be represented by a 2-coloring of the regular tiling by squares. In this paper, we will assume that the Batak weaves under consideration are periodic.

**Illustration 3.1**. To illustrate the concept of the design of a weave, let us consider the weave pattern $W_1$ appearing at the front side of the Batak basket given in Fig. 3(a), and assume it is periodic. Its design $D(W_1)$ as a 2-coloring of the regular tiling by squares is shown in Fig. 3(b). □

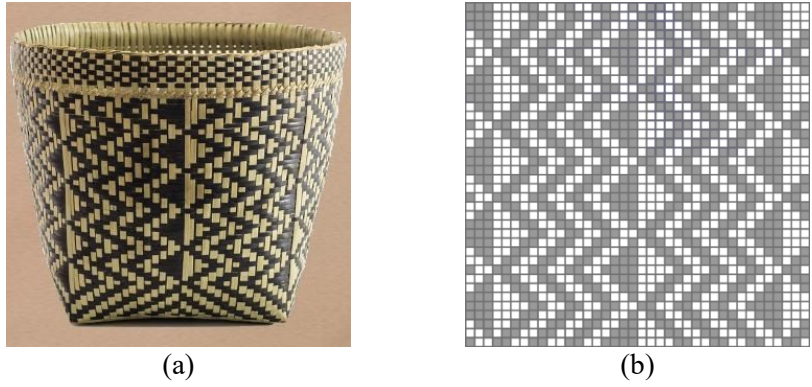

Figure 3: (a) Batak basket showing a particular weave $W_1$; and (b) the design $D(W_1)$ corresponding to $W_1$.

To determine the group of symmetries of $W$ mathematically, we will employ the setting of analyzing a 2-way 2-fold fabric given in Roth (1993) as follows. We fix $D(W)$ on a plane P, described by $z = 0$ in three-dimensional space. The weave $W$ is then considered to lie above and below P at a certain uniform distance, the warp strands "lying parallel" to the rows of $D(W)$ and the weft strands also "lying parallel" to the columns of $D(W)$.

A *symmetry* of $W$ is an isometry, a distance preserving transformation of the Euclidean space $\mathbb{E}^3$, that sends the weave $W$ to itself. The *symmetry group* $G$ of $W$ consists of all the symmetries of $W$. This is a three-dimensional group of isometries leaving the plane P invariant. $G$ is a *layer group*, which is a group of isometries of $\mathbb{E}^3$ that has translations in two directions. From Roth (1993), $G$ is embedded in a group $U = \bar{G} \times \{e, \tau\}$ where $e$ is the identity isometry of $\mathbb{E}^3$, $\tau$ is a reflection about P and $\bar{G}$ is the maximal group of symmetries of $D(W)$, that is, $\bar{G} \cong p4mm$, a plane crystallographic group. The symmetry group $G$ that leaves $W$ or $D(W)$ invariant is a subgroup of $U$. An element of $G$ is either of the form $(g, e)$ or $(g, \tau)$ where $g \in \bar{G}$. Since $\tau$ is not in $G$, an element $g$ appears in at most one pair. The set $S = \{g \in \bar{G} : g \text{ in either } (g, e) \text{ or } (g, \tau)\}$ is a subgroup of $\bar{G}$. The set $S_1 = \{g \in \bar{G} : g \text{ in } (g, e)\}$ is a subgroup of index 1 or 2 in $S$. The layer group $G$ is described by the pair $(S, S_1)$ where $S = S_1 \cup S_2$, $S_2 = \{g \in \bar{G} : g \text{ in } (g, \tau)\}$ (Roth, 1993). If $S = S_1$ the notation used is $(S, -)$. The elements of $S_1$ correspond to the elements of $G$ that preserve the side of the fabric; those in $S_2$ correspond to elements of $G$ that reverse the side of the fabric.

**Illustration 3.2** Let us assume the weave $W_1$ presented earlier lies on a plane P. We discuss some symmetries of $W_1$. $W_1$ has symmetry $(\rho_1, e)$ where $\rho_1$ is a reflection with red axis shown in Fig. 4(a). The effect of $(\rho_1, e)$ is that it preserves the side of $W_1$, that is a side of $W_1$ is sent to itself. Observe that the highlighted warp strands are sent to warp strands so the symmetry is indeed side preserving. Another symmetry of $W_1$ is $(\rho_2, \tau)$ where $\rho_2$ is a reflection with blue axis shown in Fig. 4(b) followed by a reflection $\tau$ about the plane P. The effect of $(\rho_2, \tau)$ is that it sends the front of $W_1$ to its back and reverses the side of $W_1$. Observe that the highlighted region in green consisting of warp strands is sent to the highlighted region in yellow consisting of weft strands so this symmetry is indeed side reversing. Other examples of symmetries of $W_1$: $(\sigma_1, e)$ (side preserving) where $\sigma_1$ is a glide reflection with red axis shown in Fig. 4(c), $(\sigma_2, \tau)$ (side reversing) where $\sigma_2$ is a glide reflection with blue axis given in Fig. 4(d) and $(\mu_2, \tau)$ (side reversing) where $\mu_2$ is a 2-fold rotation with blue center shown in Fig. 4(e). To illustrate other symmetries, we consider another weave $W_2$. The symmetries present in $W_2$ are: $(\mu_1, e)$ (side preserving) $\mu_1$ a 2-fold rotation with red center in Fig. 4(f), $(\lambda_1, e)$ (side preserving) where $\lambda_1$ is a translation with red vector given in Fig. 4(g) and $(\lambda_2, \tau)$ (side reversing) where $\lambda_2$ is a translation with blue vector shown in Fig. 4(h). We highlight regions of the weaves in green and yellow as a guide in Figs. 4(a)-(h). □

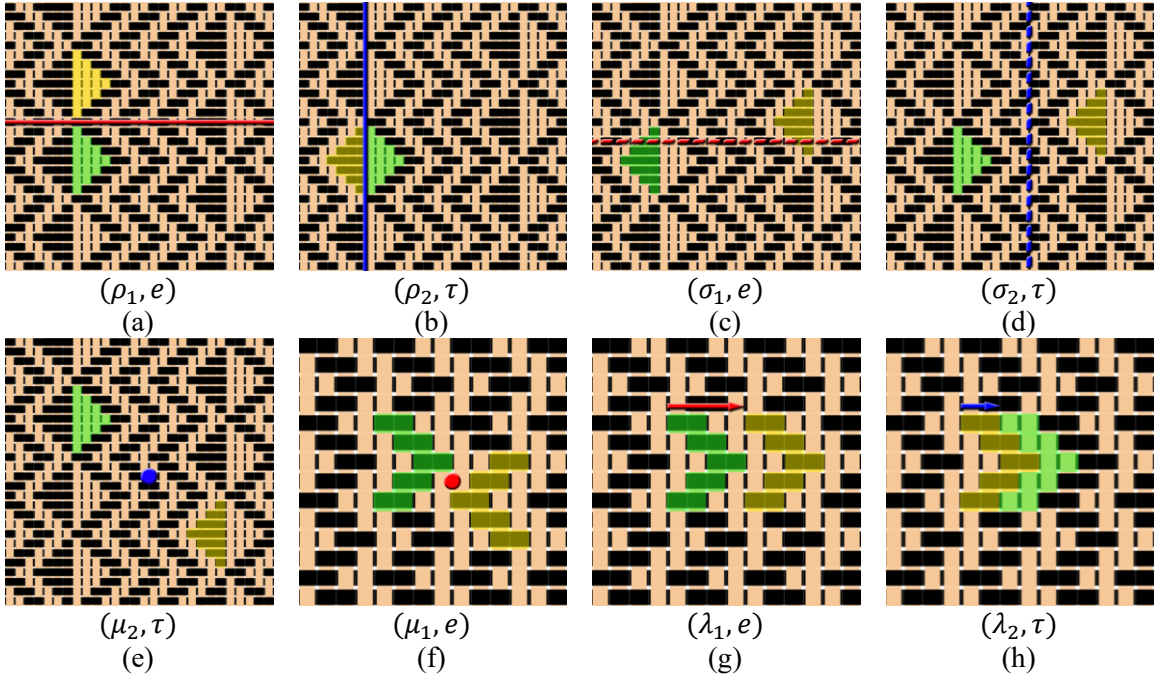

Figure 4: Different kinds of symmetries present in two examples of Batak weaves.

To describe the symmetry group $G$ of a weave $W$ using the pair $(S, S_1)$ where $S = S_1 \cup S_2$ the following result from De Las Peñas *et al.* (2023) is helpful in determining $S$.

**Theorem.** Consider a 2-way 2-fold weave $W$ with design $D(W)$. If $D(W)$ has color group $G^*$, then $G^* = S$, where $S = \{g \in \bar{G}: g$ in either $(g, e)$ or $(g, \tau)$, elements of $G\}$; $G$ is the symmetry group of $W$, $\bar{G}$ is the symmetry group of the uncolored tiling by squares.

Note that the design $D(W)$ of a weave $W$ is a 2-coloring of the regular tiling by squares so from color symmetry theory we can refer to its color group, which we define below.

Consider the set $X$ of square tiles of a regular tiling $\mathcal{T}$ by squares with symmetry group of the uncolored tiling given by $\bar{G}$. If $C = \{$grey, white$\}$, the set of two colors, a 2-coloring of $\mathcal{T}$ is a surjective function $f: \mathcal{T} \to C$ which assigns to each $x \in X$ a color $f(x)$ in $C$. Let $G^*$ be the subgroup of $\bar{G}$ that effect a permutation of the colors in $C$ called the *color group* of the given coloring of $\mathcal{T}$. Then $g^* \in G^*$ if for every $c_i \in C$ there is an element $c_j \in C$ such that $g^*(f^{-1}(c_i)) = f^{-1}(c_j)$.

In this paper, the color group $G^* = S$ of the design $D(W)$ will be illustrated geometrically by considering its diagram. The *diagram* of the color group $G^*$ is the union of the set consisting of the centers of rotations of $G^*$ and the set consisting of the axes of reflections and glide reflections of $G^*$. From the color group $G^* = S$ of $D(W)$, we determine the group $S_1$ by determining all the elements of $S$ that correspond to the side preserving symmetries of $G$.

**Illustration 3.3.** Consider the design of weave $W_1$ given in Fig. 3(b). Its color group $G_1^*$ consisting of all elements of $\bar{G}_1$ that effect a permutation of the colors is given by the plane crystallographic group $c2mm$. The diagram of $G_1^* \cong c2mm$ superimposed in $D(W_1)$ is illustrated in Fig. 5(a). From the color group $G_1^*$ we determine the index 2 subgroup consisting of elements that correspond to the elements of $S_1$. In this case we find that $S_1 \cong c1m1$. In the diagram of $G_1^*$, we color the fixed points (centers) and fixed lines (axes of symmetries) of elements of $S_1$ or $S_2$, red or blue respectively. (We will use this coloring scheme throughout the paper). For instance, $S_1$ contains the reflections and glide reflections with horizontal axes so these axes are colored red. Those axes corresponding to the reflections and glide reflections with vertical axes, as well as the centers of 2-fold rotations which are elements of $S_2$, are colored blue (Fig. 5(b)). We only show a portion of the diagram in Fig. 5(c) that corresponds to the unit cell of $S_1$. The symmetry group of $W_1$ is described by $G_1 \coloneqq (c2mm, c1m1)$. □

Now given a weave with symmetry group $G$ described by $(S, S_1)$, the next step would be to characterize its layer group structure. This is carried out by determining the three-dimensional isometries corresponding to each of the symmetries given in the diagram of the color group $G^* = S$. Then the appropriate lattice diagram of the layer group is drawn and the type of layer group can be identified from Kopský & Litvin (2002). We use Table 1 which is obtained from De Las Peñas *et al*. (2023) for our analysis based on the symmetries present in the color groups of the designs representing the Batak weaves.

Table 1: Correspondence of elements of S with three dimensional symmetries

| **Element of S** | **Three-dimensional symmetry** |
|---|---|
| Translation with vector $v$ (element of $S_1$) | Translation with vector $v$ |
| Translation with vector $v$ (element of $S_2$) | Glide reflection with glide component $v$ through the plane P |
| 2-fold rotation with center $A$ (element of $S_1$) | 2-fold rotation with axis normal to P at $A$ |
| 2-fold rotation with center $A$ (element of $S_2$) | Inversion with center $A$ |
| Reflection with axis $l$ (element of $S_1$) | Reflection with plane normal to P at $l$ |
| Reflection with axis $l$ (element of $S_2$) | 2-fold rotation with axis $l$ |
| Glide reflection with glide vector $v$ and axis $l$ (element of $S_1$) | Glide reflection with glide component $v$ through the plane normal to P at $l$ |
| Glide reflection with glide vector $v$ and axis $l$ (element of $S_2$) | 2-fold screw rotation with screw axis $v$ and axis $l$ |

**Illustration 3.4.** We now describe the symmetry group of the weave $W_1$ given earlier in Fig. 3(a) as a layer group. Let us consider the design $D(W_1)$ of the weave in Fig. 5(c) with the given diagram of its color group $S = G_1^*$. The 2-dimensional isometries present in the lattice diagram in Fig. 5(c) are analyzed with respect to their 3-dimensional counterparts as stated in Table 1. The 3-dimensional isometries are exhibited in Fig. 5(d) as follows. The reflections with red axes correspond to reflection planes perpendicular to P passing through axes represented by solid black lines. Similarly, the glide reflections with dashed red axes correspond to glide reflection planes perpendicular to P passing through axes represented by dashed black lines. The 2-fold rotations with centers at the blue circles correspond to inversions in P with centers denoted by hollow circles. The reflections with vertical blue axes correspond to 2-fold rotations in P with axes indicated by lines with arrowheads. Finally, the glide reflections correspond to 2-fold screw rotations in P with axes indicated by lines with half arrowheads. With these correspondences, we are able to construct the lattice diagram of the layer group $c2/m11$ (Fig. 5(d)) corresponding to $G_1^* \cong c2mm$ (Fig. 5(c)). □

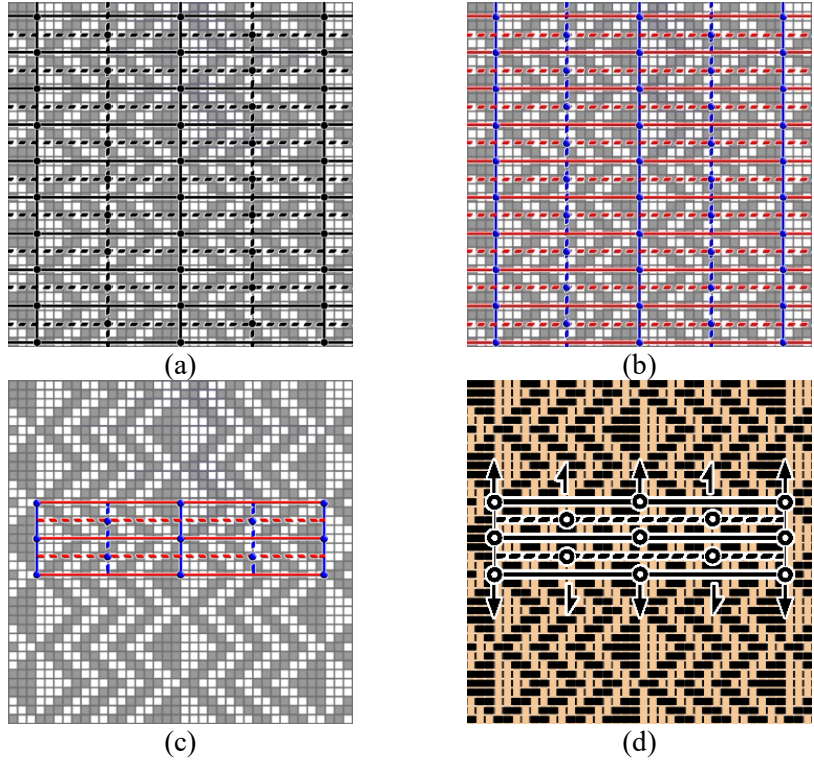

Figure 5: (a)-(b) The fixed points and fixed lines of elements of the color group $G_1^*$; (c) diagram of $G_1^* \cong c2mm$; and (d) lattice diagram of $c2/m11$.

## IV. Symmetry Groups of the Batak weaves

Twenty-two black and white motifs in Batak baskets were catalogued by Novellino where a pattern may include variations (Novellino, 2009). These patterns are based on cultural tradition and have been passed on through several generations (Novellino, 2009). Twenty-nine patterns were recorded by Calderon (1986). In our studies, we found 33 patterns in baskets and 11 in trays and mats with distinct motifs (see also De las Peñas, *et al.*, 2021). To refer to a particular type of pattern or its variant, the Batak weaver usually gives the patterns names. In our research, 22 patterns have been given names by the weavers themselves and Batak advocates.

The types of motifs and how these are repeated in the weave contribute to the symmetry structure present in the weave. In De Las Peñas *et al.* (2023) it is known that there are exactly 50 layer groups which are possible symmetry groups structures of 2-way 2-fold weaves. Our results show that 15 of these layer groups occur as symmetry groups of the Batak weaves. We have seen a prevalence of glide reflection symmetries, about 72.72% of all the patterns have symmetry groups with glide reflections.

In Fig. 6 we present the results for the basket weaves and in Fig. 7 we present the results for the non-basket weaves. The description of each symmetry group $G$ in terms of the pair $(S, S_1)$ and the layer group structure is given. Most of the geometric elements present in the weaves are triangles and diagonal lines. Diagonal lines are particularly prominent in patterns such as *Piyatwad*, *Libo-libo*, *Tiagudua*, and *Biyanig*. Additionally, some patterns incorporate both triangles and diagonal lines, such as *Enasang*, *Peitak-pietak*, *Piyakdan*, *Tiningkulob*, and *Pianpo*. In most cases triangular motifs come in pairs of black and white strands to form a diamond motif.

| | **Batak Weave W** | **Lattice Diagram of S on D(W)** | **Lattice Diagram of Layer Group on the Idealized Weave of W** |
|---|---|---|---|
| 1 | 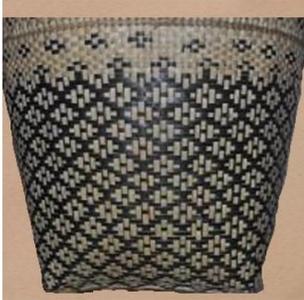 | 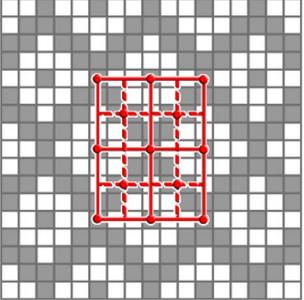 $(c2mm,-)$ | 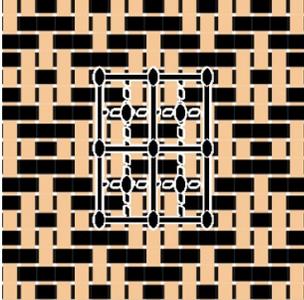 $cmm2$ |
| 2 | 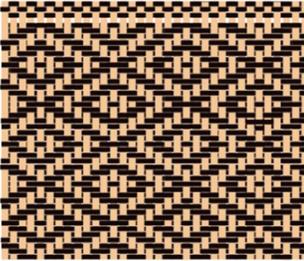 *Kiyubo* | 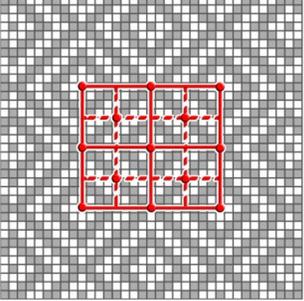 $(c2mm,-)$ | 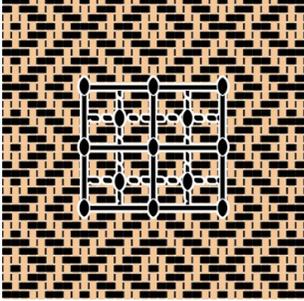 $cmm2$ |
| 3 | 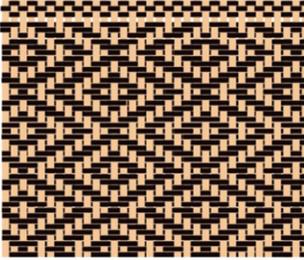 *Piyaglipusan* | 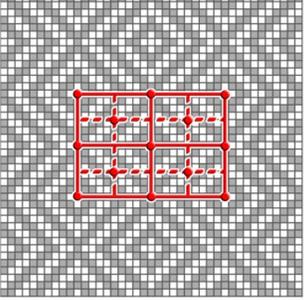 $(c2mm,-)$ | 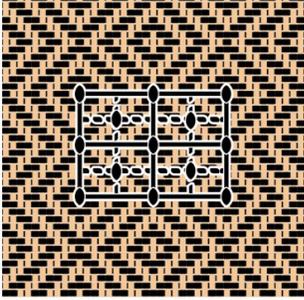 $cmm2$ |
| 4 | 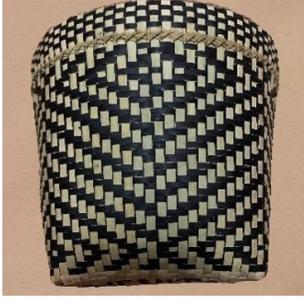 *Piyatwad/ Giyanggangan bialingan* | 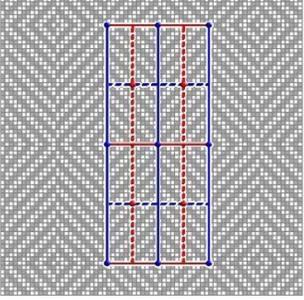 $(c2mm, p2mg)$ | 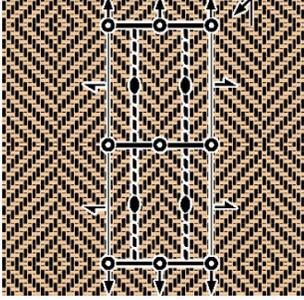 $pman$ |
| 5 | 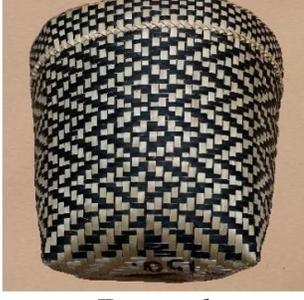 *Timograk* | 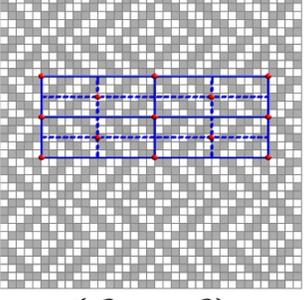 $(c2mm, p2)$ | 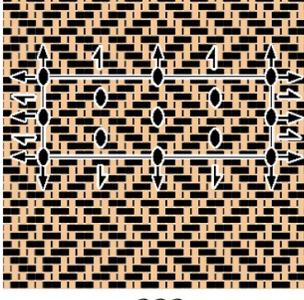 $c222$ |

| 6 | 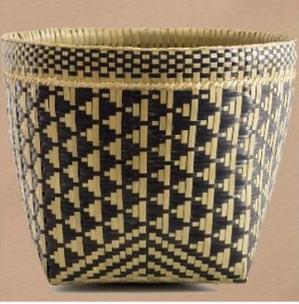 | 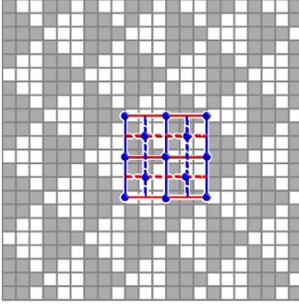 $(c2mm, c1m1)$ | 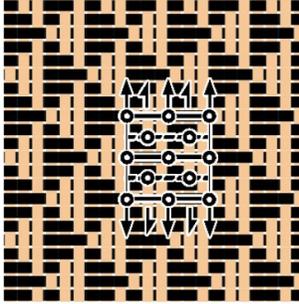 $c2/m\,11$ |
|---|---|---|---|
| 7 | 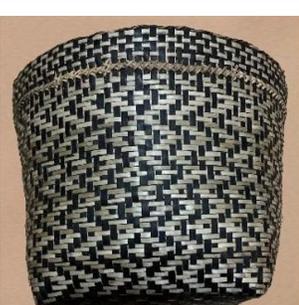 | 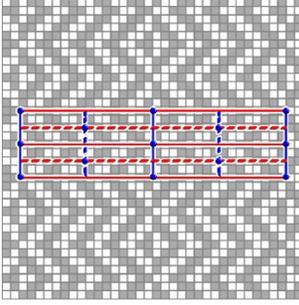 $(c2mm, c1m1)$ | 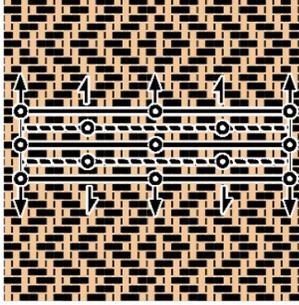 $c2/m\,11$ |
| 8 | 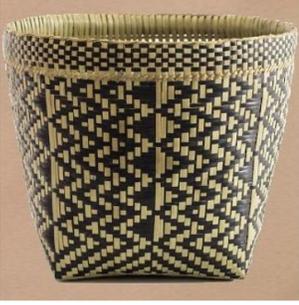 | 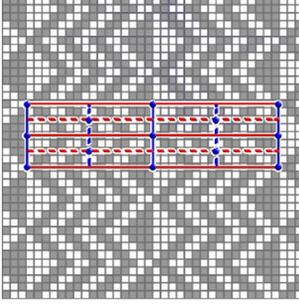 $(c2mm, c1m1)$ | 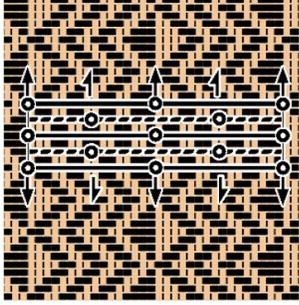 $c2/m\,11$ |
| 9 | 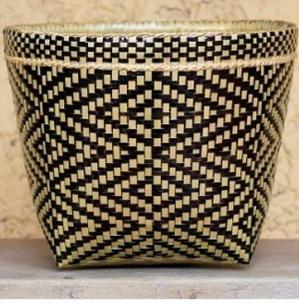 *Giyanggangan* | 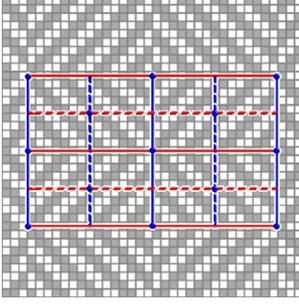 $(c2mm, c1m1)$ | 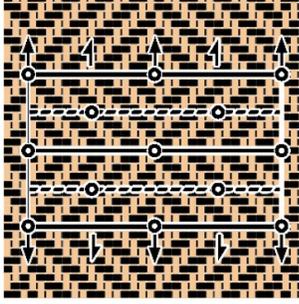 $c2/m\,11$ |
| 10 | 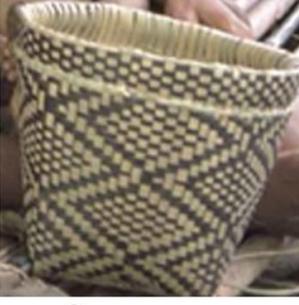 *Giyanggangan piyangapaulan* | 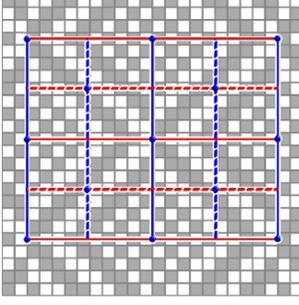 $(c2mm, c1m1)$ | 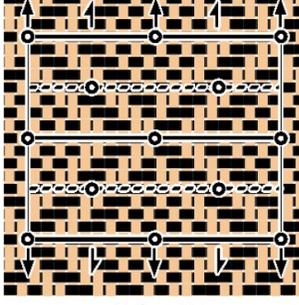 $c2/m\,11$ |
| 11 | 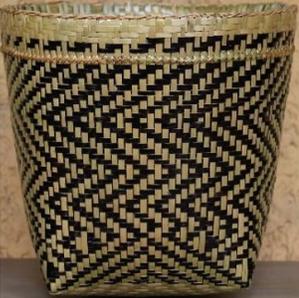 *Liangob-liangob* | 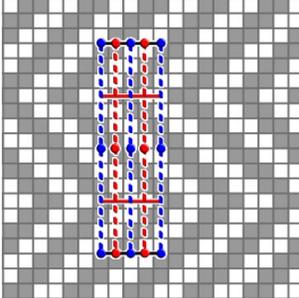 (*p2mg, p2mg*) | 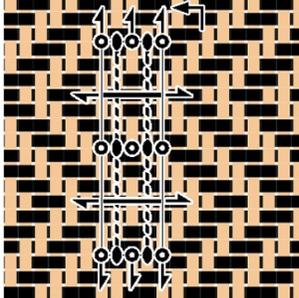 *pmab* |
| 12 | 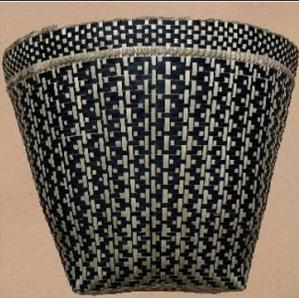 *Liolo* | 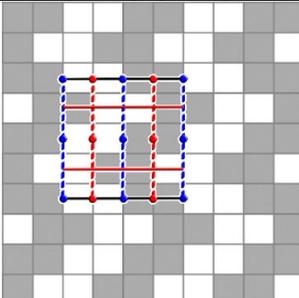 (*p2mg, p2mg*) | 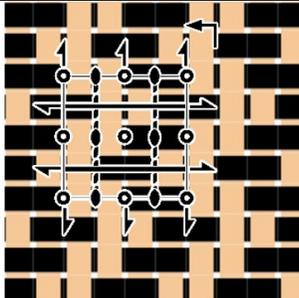 *pmab* |
| 13 | 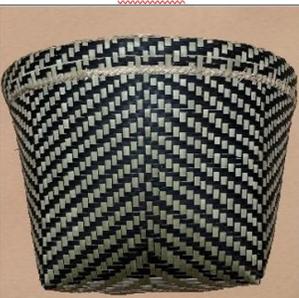 *Libo-libo* | 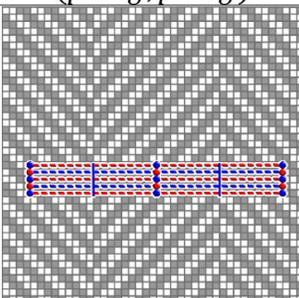 (*p2mg, p2gg*) | 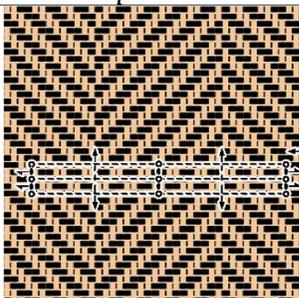 *pbab* |
| 14 | 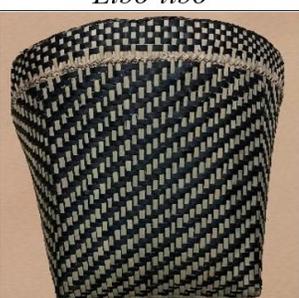 *Tiagudua* | 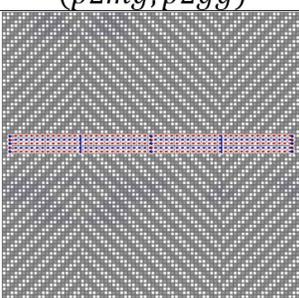 (*p2mg, p2gg*) | 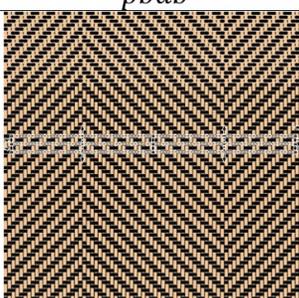 *pbab* |
| 15 | 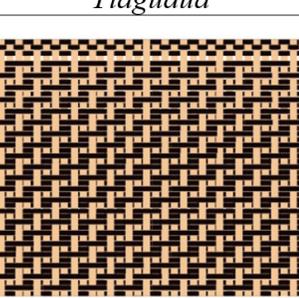 *Koyukoy* | 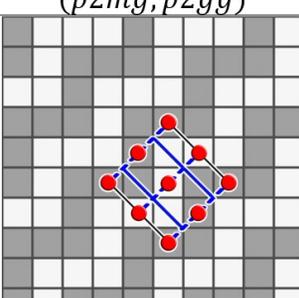 (*p2mg, p211*) | 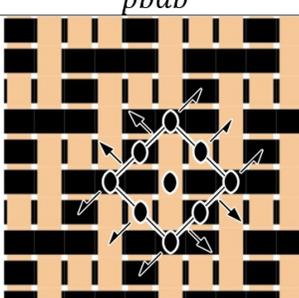 $p2_122$ |
| 16 | 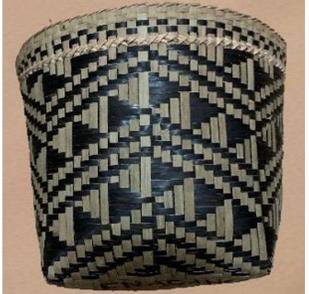<br>*Enasang* | 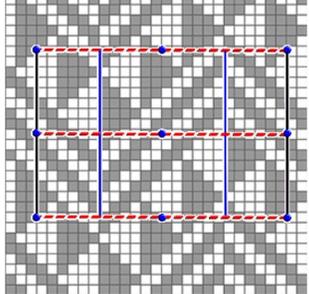<br>($p2mg, p1g1$) | 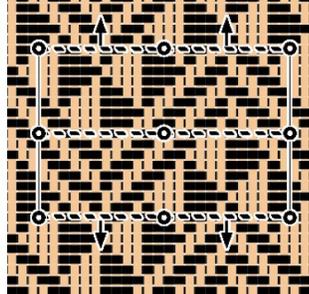<br>$p2/b\,11$ |
| 17 | 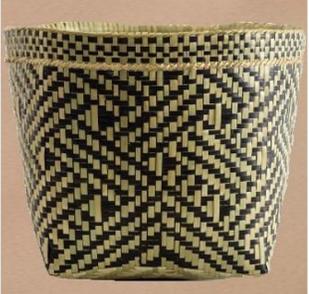 | 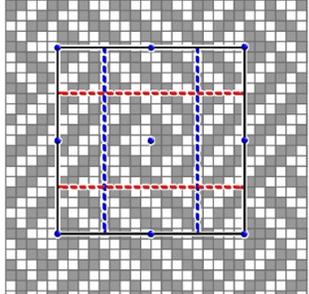<br>($p2gg, p1g1$) | 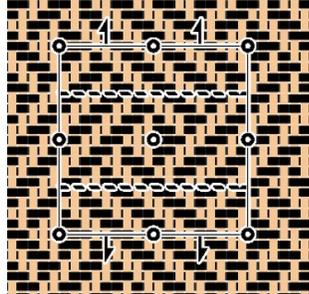<br>$p2_1/b\,11$ |
| 18 | 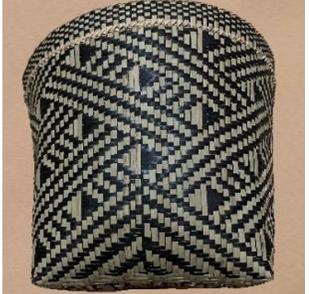<br>*Pietak pietak* | 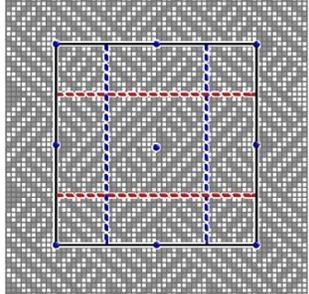<br>($p2gg, p1g1$) | 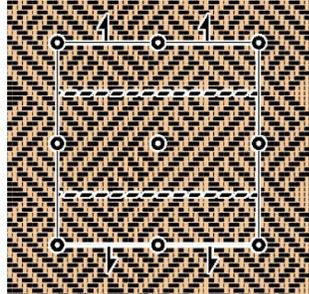<br>$p2_1/b11$ |
| 19 | 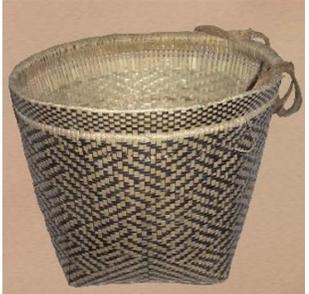<br>*Natagainpun* | 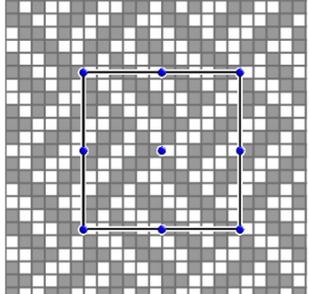<br>($p2, p1$) | 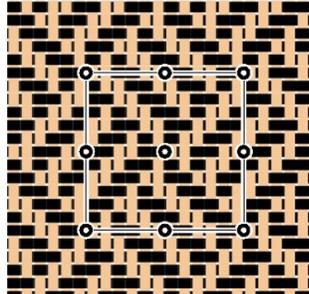<br>$p\bar{1}$ |
| 20 | 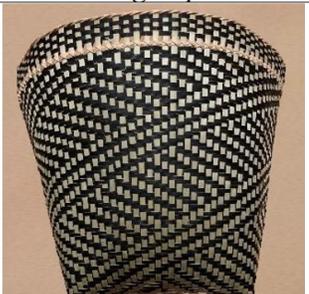<br>*Biyanig* | 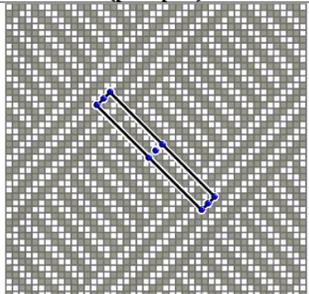<br>($p2, p1$) | 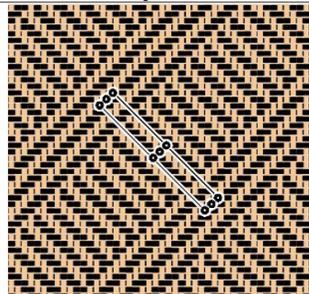<br>$p\bar{1}$ |

| 21 | 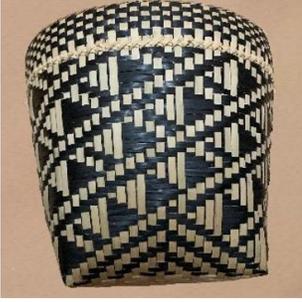 | 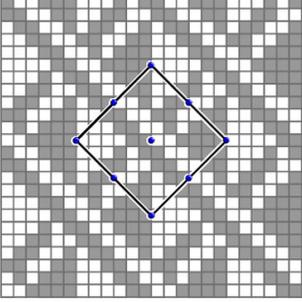 (*p*2, *p*1) | 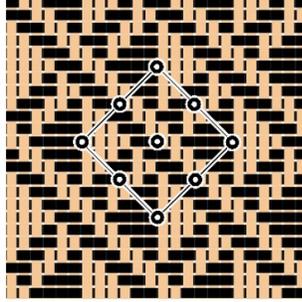 $p\bar{1}$ |
|---|---|---|---|
| 22 | 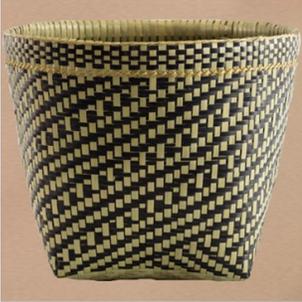 | 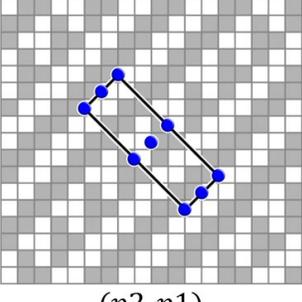 (*p*2, *p*1) | 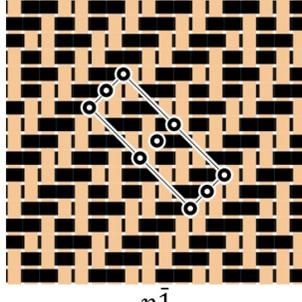 $p\bar{1}$ |
| 23 | 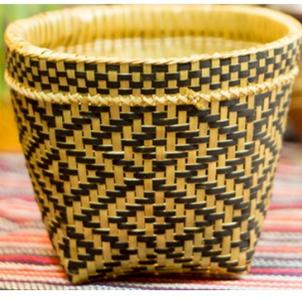 | 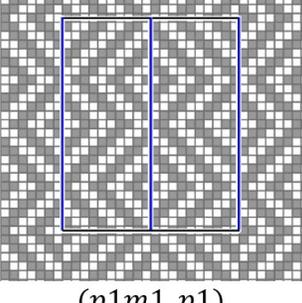 (*p*1*m*1, *p*1) | 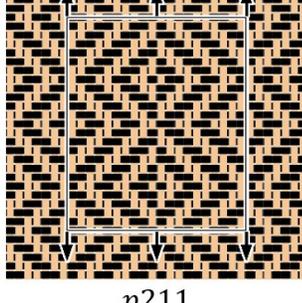 *p*211 |
| 24 | 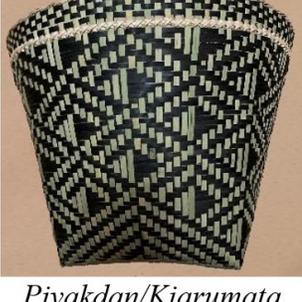 *Piyakdan/Kiarumata* | 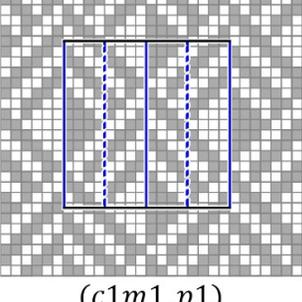 (*c*1*m*1, *p*1) | 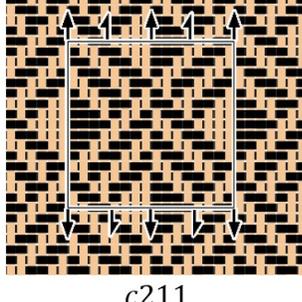 *c*211 |
| 25 | 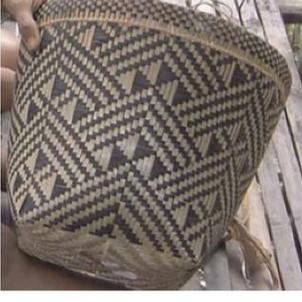 | 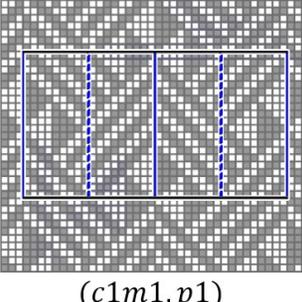 (*c*1*m*1, *p*1) | 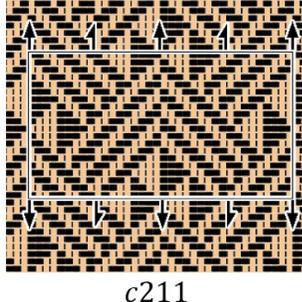 *c*211 |

| 26 | 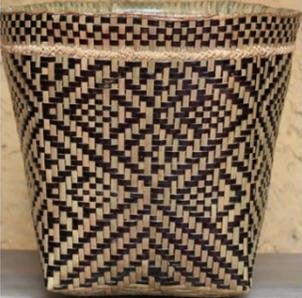 <br> *Piyaglipusan* variant | 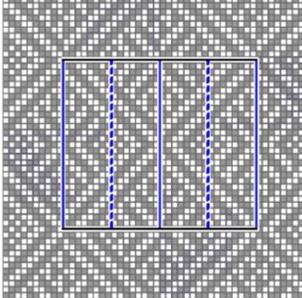 <br> $(c1m1, p1)$ | 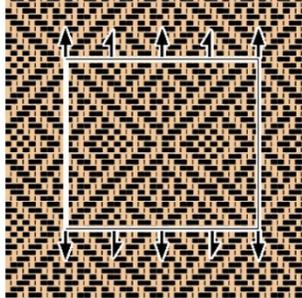 <br> $c211$ |
|---|---|---|---|
| 27 | 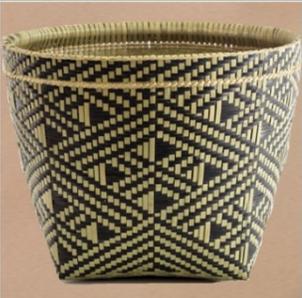 <br> *Tiningkulob* | 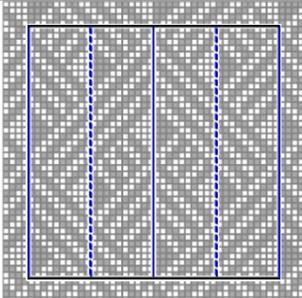 <br> $(c1m1, p1)$ | 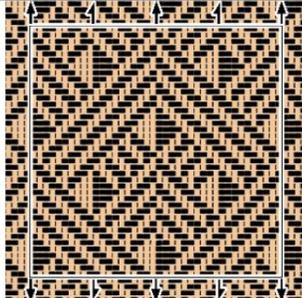 <br> $c211$ |
| 28 | 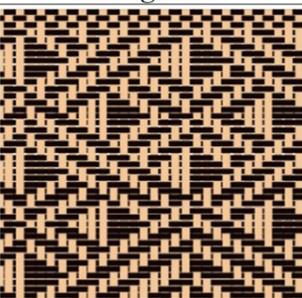 <br> *Tiningkulob* variant | 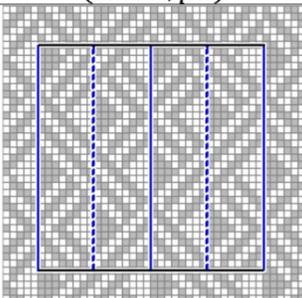 <br> $(c1m1, p1)$ | 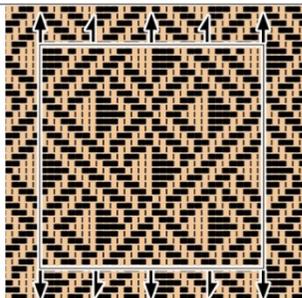 <br> $c211$ |
| 29 | 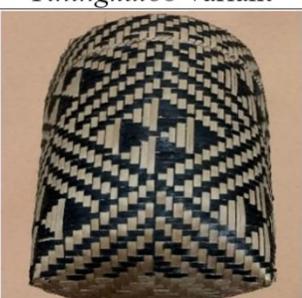 | 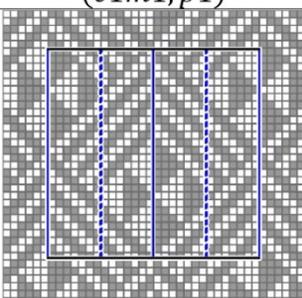 <br> $(c1m1, p1)$ | 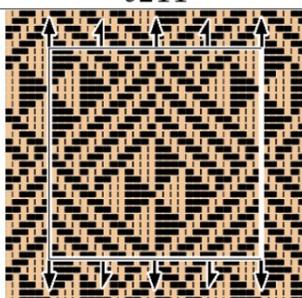 <br> $c211$ |
| 30 | 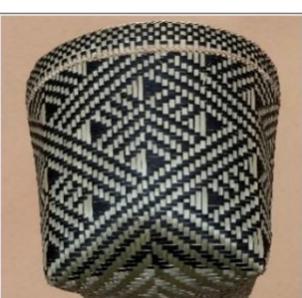 <br> *Pianpo* | 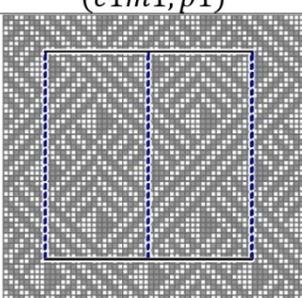 <br> $(p1g1, p1)$ | 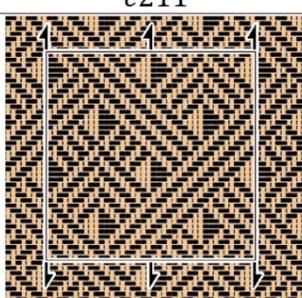 <br> $p2_111$ |

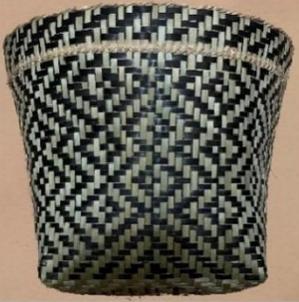

Figure 6: The symmetry groups of Batak weaves present in baskets. The first column shows either a basket exhibiting the weave at the front, or a portion of the weave; the second column gives the lattice diagram of the color group on the design of the weave with the corresponding description of $(S, S_1)$ and the third column gives the lattice diagram of the layer group in the idealized weave with the corresponding layer group structure.

| | Batak Weave $W$ | Lattice Diagram of $S$ on $D(W)$ | Lattice Diagram of Layer Group on the idealized Weave $W$ |
|---|---|---|---|
| 1 | 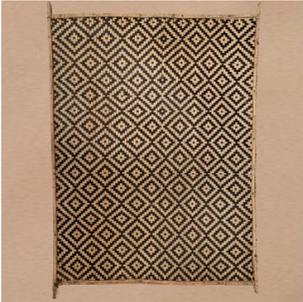 *Kiyubo* | 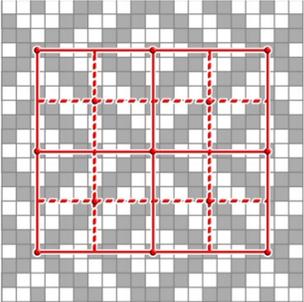 $(c2mm, -)$ | 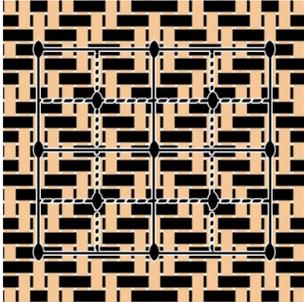 $cmm2$ |
| 2 | 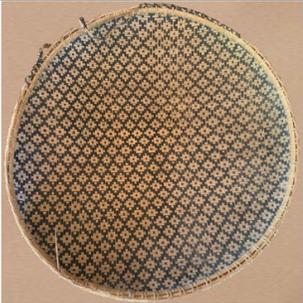 *Kiyubo* Variant | 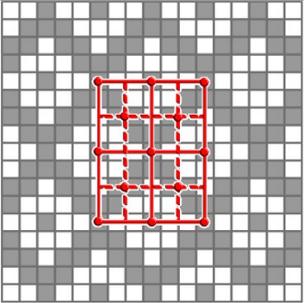 $(c2mm, -)$ | 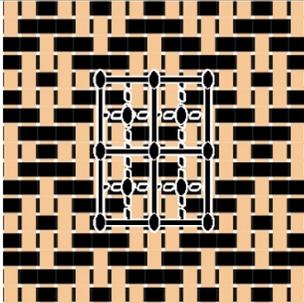 $cmm2$ |
| 3 | 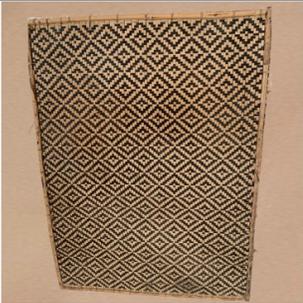 *Giyanggangan* | 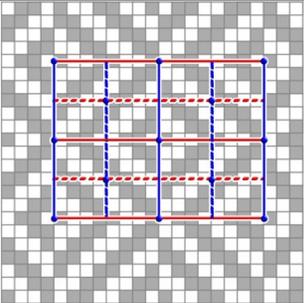 $(c2mm, c1m1)$ | 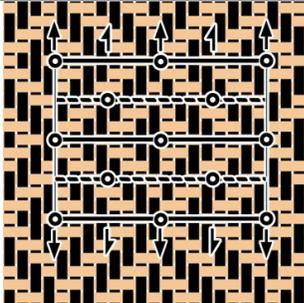 $c2/m11$ |
| 4 | 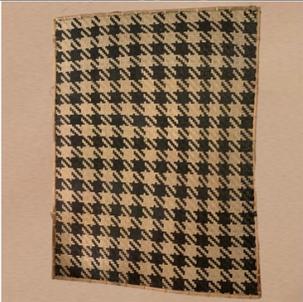 | 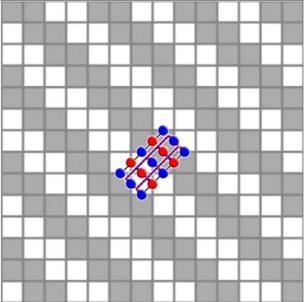 $(p2mg, p2gg)$ | 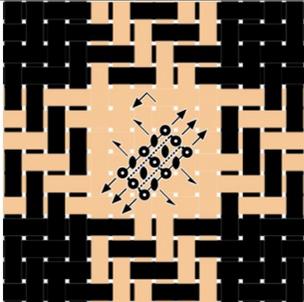 $pbab$ |
| 5 | 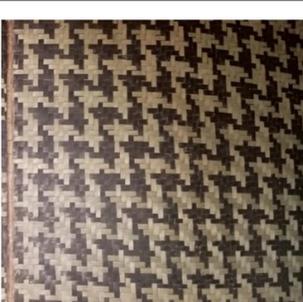 | 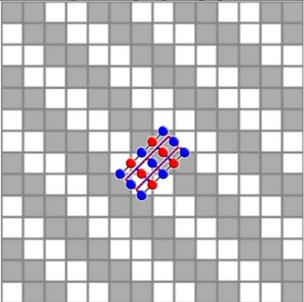 $(p2mg, p2gg)$ | 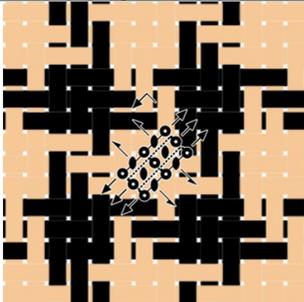 $pbab$ |
| 6 | 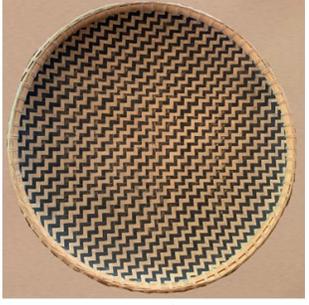 *Koyukoy* | 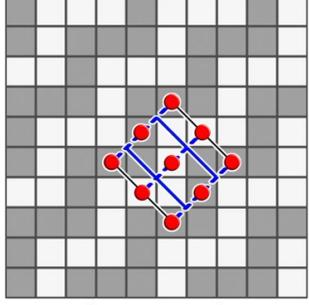 $(p2mg, p2)$ | 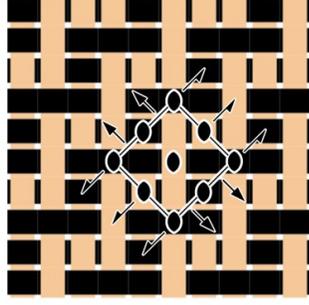 $p2_122$ |
| 7 | 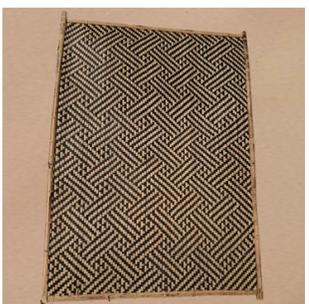 | 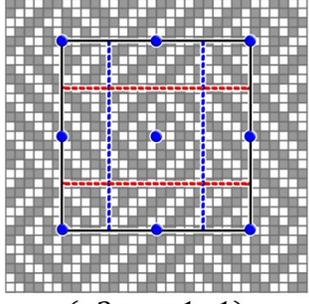 $(p2gg, p1g1)$ | 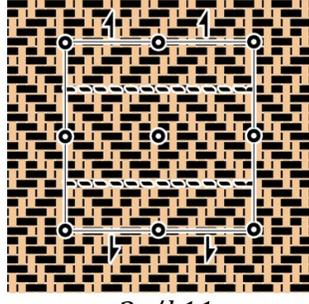 $p2_1/b11$ |
| 8 | 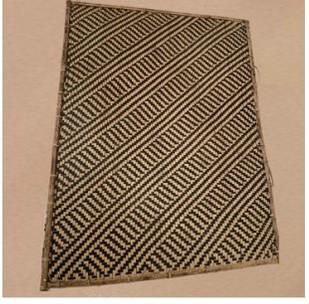 *Biyanig* | 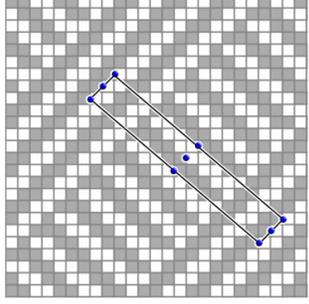 $(p2, p1)$ | 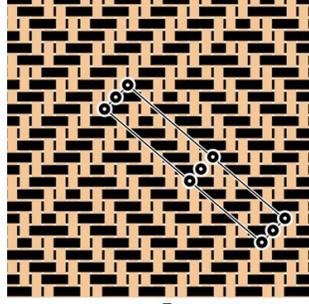 $p\bar{1}$ |
| 9 | 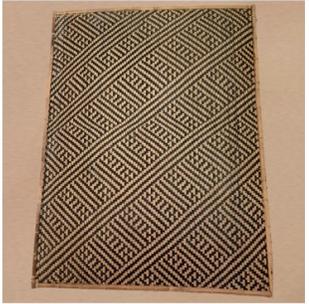 | 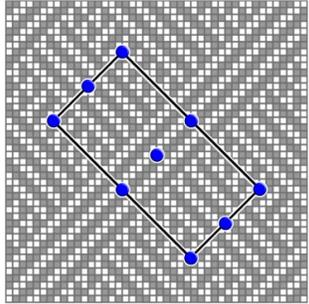 $(p2, p1)$ | 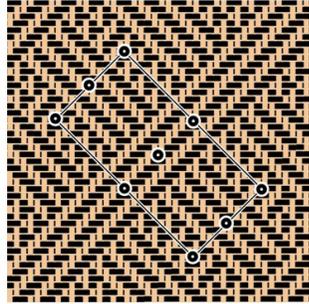 $p\bar{1}$ |
| 10 | 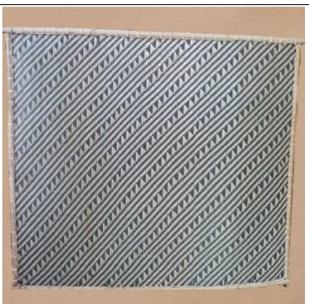 | 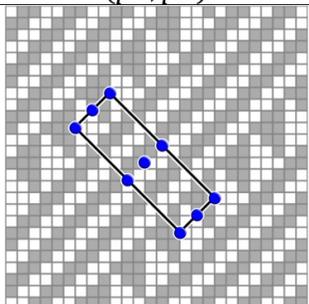 $(p2, p1)$ | 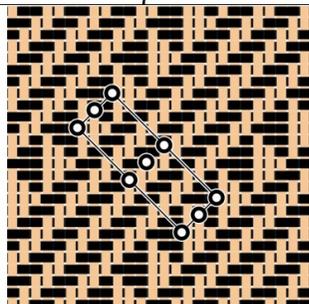 $p\bar{1}$ |

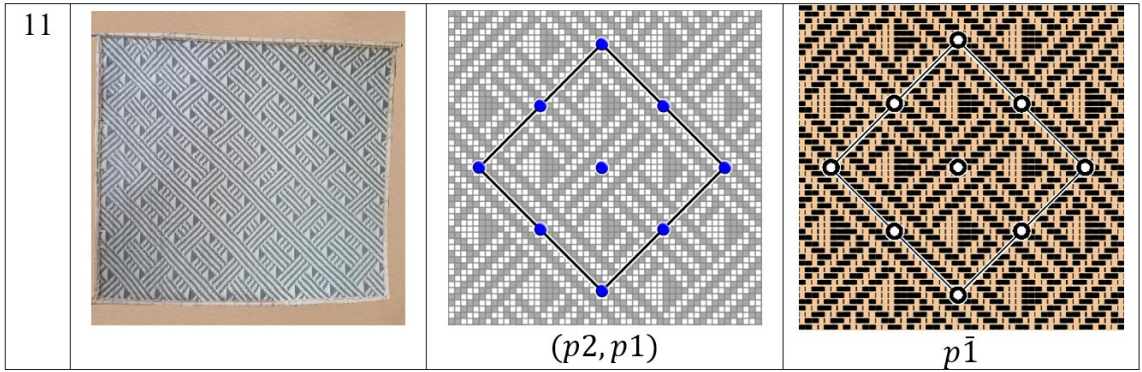

Figure 7: The symmetry groups of Batak weaves present in mats and trays. The first column shows either a mat or a tray exhibiting the weave; the second column gives the lattice diagram of the color group on the design of the weave with the corresponding description of $(S, S_1)$ and the third column gives the lattice diagram of the layer group in the idealized weave with the corresponding layer group structure.

## V. Concluding Remarks and Future Direction

This paper discusses the symmetry structures of the repeating patterns present in the baskets and other household items woven by the Batak community in Palawan. In particular, it reports the occurrence of 15 layer group structures that are present in the Batak weaves. We have observed that about 73% of the patterns studied contain glide reflection symmetries.

Several factors contribute to the occurrence of these groups. One is the use of the 2-way 2-fold weaving method and the other is due to Batak weaving techniques. The presence of these groups is also due to the types of repeating patterns present which are influenced by motifs that are culturally significant to the Batak, and have been passed on from generation to generation. For example, from the results in Figs. 6 and 7, it can be verified that the rotational symmetries present are 2-fold rotations. We have not found a 4-fold rotational symmetry in the Batak weaves even if this symmetry can occur in a 2-way 2-fold weave (see an example in Figs. 8(a)-(b)). The Batak motifs and how these are repeated in the patterns using their weaving method, do not allow for the strands to meet at right angles in the same way exhibited in Fig. 8b. The absence of the 4-fold rotational symmetries thus, does not allow for the presence of 11 layer groups as symmetry groups particularly those with four 4-rotational symmetries or 4-fold rotoinversion symmetries.

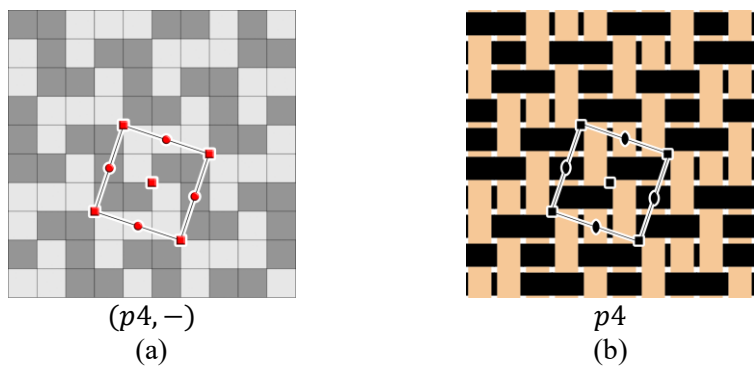

(p4, −)     p4
(a)          (b)

Figure 8: Example of a (a) design; and (b) its corresponding weave with 4-fold rotational symmetry not found in Batak weaves.

It is apparent from the results of the study that the Batak can create symmetric black and white patterns, whether this is repeated in four sides of a basket or whether these appear in the front and side of trays and mats. Very remarkable is not just the counting methods that the Batak employ to arrive at these patterns, but also the calculated technique of inserting one-sided black strips to create the patterns as a form of cultural expression. The patterns gathered in this research could be a starting point to produce a color symmetry framework in weaves using one-sided colored strands and add to existing literature on colored fabrics (e.g. Roth, 1995; Thomas, 2013).

This study shows the realization of crystallographic groups in cultural artifacts. A possible next step is to continue studies of symmetries in culture through the analysis of the weave structures that occur in other types of basket weaving traditions, including those that employ a realization of other weaves, such as the three-way three-fold weaves. It would be interesting to see how the symmetries in patterns on baskets from other cultures compare with those found in the Batak weaves.

**Acknowledgements.** Thanks go to Ms. Lara Frayre and Mr. Renato Estepa Jr. for facilitating the authors' site visit to the Batak weaving community and for sharing their pictures of the Batak baskets.